\documentclass[12pt,twoside]{article}
\usepackage{amsmath}
\usepackage{amsfonts}
\usepackage{amsthm}
\usepackage{amssymb}
\usepackage{graphicx}
\usepackage{multicol}

\begin{document}
\renewcommand{\baselinestretch}{1.5}
\title{{\normalsize
{\bf Alternative proof of the ribbonness on classical link}}}
\author{{\footnotesize Akio Kawauchi}\\ 
\date{}
{\footnotesize{\it Osaka Central Advanced Mathematical Institute, Osaka Metropolitan University }}\\
{\footnotesize{\it Sugimoto, Sumiyoshi-ku, Osaka 558-8585, Japan}}\\
{\footnotesize{\it kawauchi@omu.ac.jp}}}
\maketitle
\vspace{0.25in}
\baselineskip=15pt
\thispagestyle{empty}
\newtheorem{Theorem}{Theorem}[section]
\newtheorem{Conjecture}[Theorem]{Conjecture}
\newtheorem{Lemma}[Theorem]{Lemma}
\newtheorem{Sublemma}[Theorem]{Sublemma}
\newtheorem{Proposition}[Theorem]{Proposition}
\newtheorem{Corollary}[Theorem]{Corollary}
\newtheorem{Claim}[Theorem]{Claim}
\newtheorem{Definition}[Theorem]{Definition}
\newtheorem{Example}[Theorem]{Example}

\begin{abstract} 
Alternative proof is given for an earlier presented result that if a link in 3-space bounds a compact oriented proper surface (without closed component) in the upper half 4-space, then the link bounds a ribbon surface in the upper half 4-space which is a boundary-relative 
renewal embedding of the original surface. 
\end{abstract}

\phantom{x}



\phantom{x}

\noindent{\bf 1. Introduction}

\phantom{x}

For a set $A$ in the 3-space 
${\mathbf R}^3=\{ (x,y,z)|\, -\infty<x,y,z<+\infty\}$ 
and an interval $J\subset {\mathbf R}$, let 
\[AJ=\{ (x,y,z,t)|\, (x,y,z)\in A,\, t\in J \}.\] 
The {\it upper-half 4-space} ${\mathbf R}^4_+$ is denoted by 
${\mathbf R}^3[0,+\infty)$. 
Let $k$ be a link in the 3-space ${\mathbf R}^3$, which always bounds a compact oriented proper surface $F$ embedded smoothly in the upper-half 4-space ${\mathbf R}^4_+$, where 
${\mathbf R}^3[0]$ is canonically identified with ${\mathbf R}^3$. 
Two such surfaces $F$ and $F'$ in ${\mathbf R}^4_+$ are {\it equivalent} if there is an 
orientation-preserving diffeomorphism $f$ of ${\mathbf R}^4_+$ sending $F$ to $F'$, 
where $f$ is called an {\it equivalence}. 
For a link $k_0$ in ${\mathbf R}^3$, let ${\mathbf b}$ be a band system spanning $k_0$, namely a system of finitely many disjoint oriented bands spanning the link $k_0$ in ${\mathbf R}^3$.
The pair $(k_0, {\mathbf b})$ is called a {\it banded link}. 
The {\it surgery link} of  $(k_0, {\mathbf b})$ is the link obtained from $k_0$ by surgery along 
${\mathbf b}$. 
Assume that the surgery link of a banded link $(k_0,{\mathbf b})$ is 
a trivial link $\mbox{\boldmath $\kappa$}$ in ${\mathbf R}^3$. 
Then the band system ${\mathbf b}$ is considered as a band system 
${\mbox{\boldmath $\beta$}}$ spanning $\mbox{\boldmath $\kappa$}$. 
The pair $(\mbox{\boldmath $\kappa$}, {\mbox{\boldmath $\beta$}})$ is called 
a {\it banded loop system} with {\it loop system} $\mbox{\boldmath $\kappa$}$ and 
{\it surgery link} $k_0$. 
Throughout the paper, the surgery link $k_0$ will be 
a union $k\cup{\mathbf o}$ of a link $k$ in question and a trivial link ${\mathbf o}$ called an {\it extra trivial link}. 
Here, it is assumed that there is a band sub-system ${\mathbf b}_1$ of the 
band system ${\mathbf b}$ such that ${\mathbf b}_1$ 
connects to  ${\mathbf o}$ with just one band $b_1\in{\mathbf b}_1$ for every component  $o\in{\mathbf o}$ and every band 
$b\in{\mathbf b}_1^c={\mathbf b}\setminus{\mathbf b}_1 $ spans the link $k$.
Let ${\mbox{\boldmath $\alpha$}}_1$ be the arc system of the attaching arc $\alpha_1$ 
of every band $b_1\in{\mathbf b}_1$ to $o\in{\mathbf o}$, and 
${\mbox{\boldmath $\alpha$}}_1^c$ the complementary arc system of 
${\mbox{\boldmath $\alpha$}}_1$ in ${\mathbf o}$ consisting of every 
complementary arc $\alpha_1^c=\mbox{cl}(o\setminus \alpha_1)$. 
Any disk system ${\mathbf d}$ in ${\mathbf R}^3$ bounded by the extra trivial link ${\mathbf o}$
is called an {\it extra disk system}, which is fixed and the argument proceeds.
Let $\mbox{\boldmath $\delta$}$ be a disk system 
consisting of disjoint disks in ${\mathbf R}^3$
with $\partial\mbox{\boldmath $\delta$}=\mbox{\boldmath $\kappa$}$, 
which is called a {\it based disk system} for a loop system $\mbox{\boldmath $\kappa$}$.
A ribbon surface-link $\mbox{cl}(F^1_{-1})$ in ${\mathbf R}^4$ is constructed from 
a banded loop system
$(\mbox{\boldmath $\kappa$}, {\mbox{\boldmath $\beta$}})$
by taking the surgery of the trivial $S^2$-link 
\[O=\partial (\mbox{\boldmath $\delta$}[-1,1])
=\mbox{\boldmath $\delta$}[-1]\cup (\partial\mbox{\boldmath $\delta$})[-1,1]
\cup \mbox{\boldmath $\delta$}[1]\] 
along the 1-handle system ${\mbox{\boldmath $\beta$}}[-t,t]$ in ${\mathbf R}^4$ 
for any $t$ with $0<t<1$. 
The proper surface $\mbox{ucl}(F_0^1)= \mbox{cl}(F^1_{-1})\cap{\mathbf R}^4_+$ 
in ${\mathbf R}^4_+$ is called the {\it upper-closed realizing surface} of  
a banded loop system $(\mbox{\boldmath $\kappa$}, {\mbox{\boldmath $\beta$}})$ 
with surgery link $k_0$. 
Note that choices of the based disk systems $\mbox{\boldmath $\delta$}$ are independent of 
the equivalences of $\mbox{ucl}(F_0^1)$ and $\mbox{cl}(F^1_{-1})$ by Horibe-Yanagawa's lemma, \cite{KSS}. The reason for dealing with a banded loop system 
$(\mbox{\boldmath $\kappa$}, {\mbox{\boldmath $\beta$}})$
 rather than a banded link $(k_0, {\mathbf b})$ is because not only can a based disk system 
 $\mbox{\boldmath $\delta$}$ be chosen freely, but it also makes a band deformation of 
the  band system ${\mbox{\boldmath $\beta$}}$ easier.  Actually, an isotopic deformation of  
 ${\mbox{\boldmath $\beta$}}$ respecting  the arc system ${\mbox{\boldmath $\alpha$}}_1$
 and the loop system $\mbox{\boldmath $\kappa$}$ does not change 
 the ribbon surface-link $\mbox{cl}(F^1_{-1})$ in  ${\mathbf R}^4$ and 
 the  proper surface $\mbox{ucl}(F_0^1)$ in  ${\mathbf R}^4_+$, up to equivalences.  
 
Let $\mbox{cl}(F^1_{-1})_{\mathbf d}$ be the surface-link in ${\mathbf R}^4$ 
obtained from the ribbon surface-link $\mbox{cl}(F^1_{-1})$ by surgery along 
the 2-handle ${\mathbf d}[-\varepsilon,\varepsilon]$ on $\mbox{cl}(F^1_{-1})$ 
where $0<\varepsilon<t<1$. 
The proper surface $P(F_0^1)=\mbox{cl}(F^1_{-1})_{\mathbf d}\cap {\mathbf R}^4_+$ 
in ${\mathbf R}^4_+$ with $\partial P(F_0^1)=k$ is called 
a {\it proper realizing surface} of a banded loop system
$(\mbox{\boldmath $\kappa$}, {\mbox{\boldmath $\beta$}})$ with surgery link 
$k_0=k\cup {\mathbf o}$. 
The following theorem is known, \cite{KSS}.

\phantom{x}

\noindent{\bf Normal form theorem.} Every compact oriented proper surface $F$ without closed component in the upper-half 4-space ${\mathbf R}^4_+$ with $\partial F=k$ in ${\mathbf R}^3$ 
is equivalent to a proper realizing surface $P(F_0^1)$ in ${\mathbf R}^4_+$ with 
$\partial P(F_0^1)=k$ of a banded loop system
$(\mbox{\boldmath $\kappa$}, {\mbox{\boldmath $\beta$}})$ with surgery link 
$k_0=k+{\mathbf o}$ which is a split sum of $k$ and an extra trivial link ${\mathbf o}$. 

\phantom{x}

The proper realizing surface $P(F^1_0)$ in ${\mathbf R}^4_+$ is called a 
{\it normal form} of the proper surface $F$ in ${\mathbf R}^4_+$. If the extra trivial link ${\mathbf o}$ is taken the empty link, namely 
$P(F^1_0)=\mbox{ucl}(F_0^1)$, then the proper surface $F$ in ${\mathbf R}^4_+$ 
is called a {\it ribbon surface}. 
In the following example, it is observed that there are lots of compact oriented proper surfaces without closed component in ${\mathbf R}^4_+$ which are not equivalent to any ribbon surface 
in ${\mathbf R}^4_+$. 

\phantom{x}

\noindent{\bf Example.} 
For every link $k$ in ${\mathbf R}^3$, let $F'$ be any ribbon surface in ${\mathbf R}^4_+$ 
with $k=\partial F'$. 
For example, let $F'$ be a proper surface in ${\mathbf R}^4_+$ obtained from 
a Seifert surface for $k$ in ${\mathbf R}^3$ by an interior push into ${\mathbf R}^4_+$. 
Take a connected sum $F=F'\# K$ of $F'$ and a non-trivial $S^2$-knot $K$ in ${\mathbf R}^4$ 
with non-abelian fundamental group. Then $k=\partial F'=\partial F$. 
It is shown that $F$ is not equivalent to any ribbon surface in ${\mathbf R}^4_+$. 
The fundamental groups of $k, F', F, K$ 
are denoted as follows.
\[\pi(k)=\pi_1({\mathbf R}^3\setminus k,x_0), 
\quad \pi(F')=\pi_1({\mathbf R}^4\setminus F',x_0), \]
\[\pi(F)=\pi_1({\mathbf R}^4\setminus F,x_0), \quad
\pi(K)=\pi_1(S^4\setminus K, x_0).\]
Let $\pi(k)^*, \pi(F')^*\pi(F)^*, \pi(K)^*$ be the kernels of the canonical epimorphisms 
from the groups $\pi(k), \pi(F'), \pi(F), \pi(K)$ to the infinite cyclic group sending every meridian element to the generator, respectively. 
It is a special feature of a ribbon surface $F'$ that
the canonical homomorphism $\pi(k)\to \pi(F')$ is an epimorphism, 
so that the induced homomorphism $\pi(k)^*\to \pi(F')^*$ is onto. 
On the other hand, the canonical homomorphism $\pi(k)\to \pi(F)$ 
is not onto, because the group $\pi(F)^*$ is the free product $\pi(F')^* * \pi(K)^*$
and $\pi(K)^*\ne 0$ and the image of the induced homomorphism $\pi(k)^*\to \pi(F)^*$ 
is just the free product summand $\pi(F')^*$. 
Thus, the proper surface $F$ in ${\mathbf R}^4_+$ is not equivalent to any ribbon surface. 

\phantom{x}

A compact oriented proper surface $F'$ in ${\mathbf R}^4_+$ is a {\it renewal embedding} of a compact oriented proper surface $F$ in ${\mathbf R}^4_+$ if there is an 
orientation-preserving surface-diffeomorphism $F'\to F$ keeping the boundary fixed. 
A renewal embedding $F'$ of $F$ is {\it boundary-relative} if 
the link $k'=\partial F'$ in ${\mathbf R}^3$ is equivalent to the link $k=\partial F$ in 
${\mathbf R}^3$. 
The proof of the following theorem is given, \cite{Ka2}. 
In this paper, an alternative proof of this theorem is given from a viewpoint of deformations of a ribbon surface-link in ${\mathbf R}^4$. 

\phantom{x}

\noindent{\bf Classical ribbon theorem.} Assume that a link $k$ in the 3-space 
${\mathbf R}^3$ bounds a compact oriented proper surface $F$ without closed component 
in the upper-half 4-space ${\mathbf R}^4_+$. Then the link $k$ in ${\mathbf R}^3$ 
bounds a ribbon surface $F'$ in ${\mathbf R}^4_+$ which is a boundary-relative renewal embedding of $F$.

\phantom{x}

A link $k$ in ${\mathbf R}^3$ is a {\it slice link in the strong sense} if $k$ 
bounds a proper disk system embedded smoothly in ${\mathbf R}^4_+$. 
A link $k$ in ${\mathbf R}^3$ is a {\it ribbon link} if $k$ bounds a ribbon disk system 
in ${\mathbf R}^4_+$. 
The following corollary is a special case of Classical ribbon theorem.

\phantom{x}

\noindent{\bf Corollary~1.} Every slice link in the strong sense 
in ${\mathbf R}^3$ is a ribbon link.

\phantom{x}

Thus, Classical ribbon theorem solves {\it Slice-Ribbon Problem}, \cite{Fox1}, \cite{Fox2}.
The following corollary is obtained from Corollary~1.

\phantom{x}

\noindent{\bf Corollary~2.} A link $k$ in ${\mathbf R}^3$ is a ribbon link if 
a ribbon link is obtained from the split sum $k+ {\mathbf o}$ of $k$ and a trivial link 
${\mathbf o}$ by a band sum of $k$ and every component of ${\mathbf o}$. 

\phantom{x}

The proof of the classical ribbon theorem is done throughout the section 2. 
An idea of the proof is to consider the 2-handle pair system $(D\times I,D'\times I)$ 
on the ribbon surface-link $\mbox{cl}(F^1_{-1})$ with $k+{\mathbf o}$ as the middle-cross sectional link such that $P(F^1_0)$ is 
equivalent to a previously given surface $F$ in ${\mathbf R}^4_+$, where 
the 2-handle system $D\times I$ is constructed from the band system 
${\mathbf b}_1$ and the 2-handle system $D'\times I$ is constructed from 
the extra disk system ${\mathbf d}$. 
The interior intersections of $(D\times I,D'\times I)$
will be eliminated and $(D\times I,D'\times I)$ becomes an O2-handle pair system on a new ribbon surface-link $\mbox{cl}(F^1_{-1})$ with $k+{\mathbf o}$ as the middle-cross sectional link obtained by sacrificing equivalences. 
Then $P(F^1_0)$ is a ribbon surface that is a boundary-relative renewal embedding of $F$, 
which will complete the proof.

\phantom{x}

\noindent{\bf 2. Proof of Classical ribbon theorem}

\phantom{x}

{\it Throughout this section, the proof of the classical ribbon theorem is done.} 
Let $F$ be a compact oriented proper surface without closed component 
in ${\mathbf R}^4_+$, and $\partial F= k$ a link in ${\mathbf R}^3$. 
By the normal form theorem, there is a banded loop system
$(\mbox{\boldmath $\kappa$}, {\mbox{\boldmath $\beta$}})$ with surgery link 
$k_0=k+{\mathbf o}$ such that $P(F^1_0)$ is equivalent to $F$. 
The extra trivial link ${\mathbf o}$ is uniquely specified by the banded loop system 
$(\mbox{\boldmath $\kappa$}, {\mbox{\boldmath $\beta$}})$, which is the union of 
the arc system 
${\mbox{\boldmath $\alpha$}}_1$ and the complementary arc system ${\mbox{\boldmath $\alpha$}}_1^c$, where the interior of ${\mbox{\boldmath $\alpha$}}_1$ transversely meets the interior of a based disk system $\mbox{\boldmath $\delta$}$ with finite points and is disjoint from the based loop system $\mbox{\boldmath $\kappa$}$ and 
${\mbox{\boldmath $\alpha$}}_1^c$ belongs to the loop system $\mbox{\boldmath$\kappa$}$. 

A {\it renewal embedding} of a banded loop system 
$(\mbox{\boldmath $\kappa$}, {\mbox{\boldmath $\beta$}})$ with surgery link 
$k_0=k\cup{\mathbf o}$ is a banded loop system 
$(\mbox{\boldmath $\kappa$}', {\mbox{\boldmath $\beta$}}')$ with surgery link 
$k'_0=k'\cup{\mathbf o}$  
such that  there is a  homeomorphism  
$\mbox{\boldmath $\kappa$}\cup  {\mbox{\boldmath $\beta$}}\to  
\mbox{\boldmath $\kappa$}'\cup  {\mbox{\boldmath $\beta$}}'$
with restrictios 
$\mbox{\boldmath $\kappa$}\to  \mbox{\boldmath $\kappa$}'$
and $ {\mbox{\boldmath $\beta$}}\to   {\mbox{\boldmath $\beta$}}'$ 
orientation-preserved.

\phantom{x}

The following observation is directly obtained by definition.

\phantom{x}

\noindent{\bf (2.1)} If a banded loop system 
$(\mbox{\boldmath $\kappa$}', {\mbox{\boldmath $\beta$}}')$ 
with surgery link $k'\cup{\mathbf o}$ is 
a renewal embedding of a banded loop system 
$(\mbox{\boldmath $\kappa$}, {\mbox{\boldmath $\beta$}})$ with surgery link 
$k\cup{\mathbf o}$, then the upper-closed realizing surface $\mbox{ucl}(F_0^1)'$ 
constructed from  
$(\mbox{\boldmath $\kappa$}', {\mbox{\boldmath $\beta$}}')$ is a renewal embedding of the upper-closed realizing surface $\mbox{ucl}(F_0^1)$ constructed from  
$(\mbox{\boldmath $\kappa$}, {\mbox{\boldmath $\beta$}})$ such that 
$\partial\mbox{ucl}(F_0^1)=k\cup{\mathbf o}$ and $\partial\mbox{ucl}(F_0^1)'=k'\cup{\mathbf o}$. 

\phantom{x}

A {\it transversal arc} of a band spanning a link is a simple proper arc in the band 
which is parallel to an attaching arc.
For a band $b\in{\mathbf b}$ transversely meeting the interior of an extra disk 
$d\in{\mathbf d}$, the $d$-{\it arc system} of $b$ 
is the arc system $d(b)$ of every transversal arc $a$ of $b$ in the interior of $d$. 
The ${\mathbf d}$-arc system of a band system ${\mathbf b}$ is 
the collection ${\mathbf d}({\mathbf b})$ of $d(b)$ for every $d\in{\mathbf d}$ and 
every $b\in{\mathbf b}$.
For a based disk $\delta\in\mbox{\boldmath $\delta$}$, 
the $\delta$-{\it arc system} of a band $\beta\in\mbox{\boldmath $\beta$}$ 
is the arc system $\delta(\beta)$ of every transversal arc $c$ of 
$\beta$ in the interior of $\delta$. 
The $\mbox{\boldmath $\delta$}$-arc system of $\mbox{\boldmath $\beta$}$ is 
the collection $\mbox{\boldmath $\delta$}({\mbox{\boldmath $\beta$}})$ of $\delta(\beta)$ for every $\delta\in\mbox{\boldmath $\delta$}$ and every $\beta\in\mbox{\boldmath $\beta$}$.
A {\it normal proper arc} in the extra disk system ${\mathbf d}$ is a simple proper arc in 
${\mathbf d}$ with the endpoints in the interior of the arc system 
${\mbox{\boldmath $\alpha$}}_1$. 
The following assertion is shown. 

\phantom{x}

\noindent{\bf (2.2)} By isotopic deformations in ${\mathbf R}^3$, the banded loop system 
$(\mbox{\boldmath $\kappa$}, {\mbox{\boldmath $\beta$}})$ in 
${\mathbf R}^3$ with surgery link $k_0=k+{\mathbf o}$ is deformed so that
a based disk system ${\mbox{\boldmath $\delta$}}$ transversely 
meets the extra disk system ${\mathbf d}$ with interior simple arcs or normal proper arcs in ${\mathbf d}$ except for the complementary arc system 
${\mbox{\boldmath $\alpha$}}_1^c$. 

\phantom{x}

\begin{figure}[hbtp]
\begin{center}
\includegraphics[width=7cm, height=8cm]{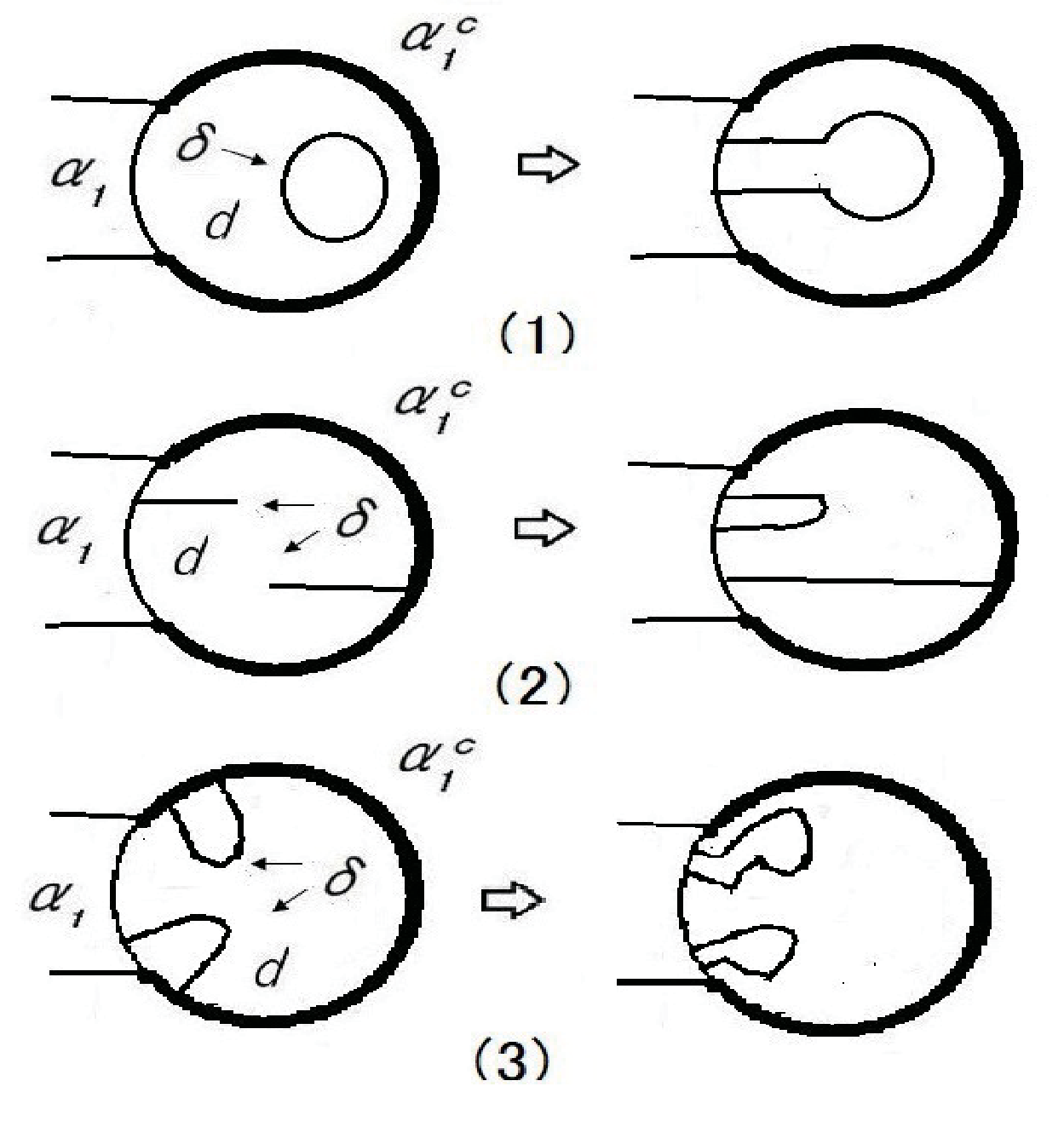}
\end{center}
\caption{Changing the intersection of a based disk and an extra disk }
\label{fig:Change}
\end{figure}

\noindent{\it Proof of (2.2).} 
By transverse regularity, 
the intersection $d\cap \delta$ for  every $d\in{\mathbf d}$ and every 
$\delta\in{\mbox{\boldmath $\delta$}}$ 
is made interior simple loops, interior simple arcs, clasp type simple arcs or 
simple proper arcs in ${\mathbf d}$ except for the complementary arc system 
${\mbox{\boldmath $\alpha$}}_1^c$. 
A simple loop is changed into a normal proper arc by a pushing out deformation 
to ${\mbox{\boldmath $\alpha$}}_1$, Fig.~\ref{fig:Change} (1). 
A clasp type simple arc is changed into a simple proper arc 
by moving out the interior point 
to ${\mbox{\boldmath $\alpha$}}_1$, Fig.~\ref{fig:Change} (2). 
A simple proper arc which is not normal is also changed into a normal proper arc
by a pushing out deformation of the arc system of ${\mbox{\boldmath $\delta$}}$ 
meeting a boundary collar of 
${\mbox{\boldmath $\alpha$}}_1^c$ in 
${\mathbf d}$, Fig.~\ref{fig:Change} (3). 
Thus, a deformed based disk system  ${\mbox{\boldmath $\delta$}}$ transversely 
meets ${\mathbf d}$ with interior simple arcs or normal proper arcs in ${\mathbf d}$  except for the complementary arc system 
${\mbox{\boldmath $\alpha$}}_1^c$.
This completes the proof of (2.2). 

\phantom{x}

The following operation gives a standard renewal embedding of a banded loop system. 

\phantom{x}

\begin{figure}[hbtp]
\begin{center}
\includegraphics[width=11cm, height=3.5cm]{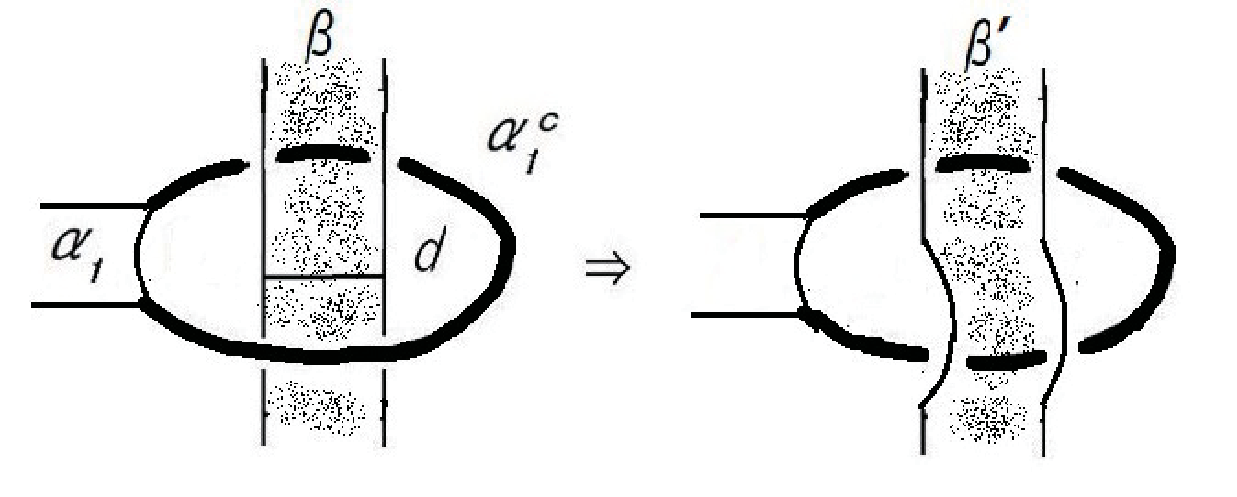}
\end{center}
\caption{Band Move Operation }
\label{fig:BMO}
\end{figure}

\noindent{\bf Band Move Operation.} 
In the banded loop system 
$(\mbox{\boldmath $\kappa$}, {\mbox{\boldmath $\beta$}})$ with surgery link 
$k_0=k\cup{\mathbf o}$, 
assume that there is a transversal arc $c$ of a 
band $\beta\in {\mbox{\boldmath $\beta$}}$ in the interior of an extra disk 
$d\in {\mathbf d}$ and 
there is a simple path $\omega$ in $d$ from a point $p\in c$ to an interior point 
of the arc $\alpha^c_1=\partial d\cap{\mbox{\boldmath $\alpha$}}_1^c$ 
which avoids meeting ${\mbox{\boldmath $\beta$}}$ other than $c$. 
Let $\beta'$ be a band obtained from $\beta$ by sliding the arc $c$ off 
the disk $d$ along the path $\omega$. Replace the banded loop system 
$(\mbox{\boldmath $\kappa$}, {\mbox{\boldmath $\beta$}})$ 
with the banded loop system 
$(\mbox{\boldmath $\kappa$}, {\mbox{\boldmath $\beta$}}')$ 
obtained by replacing $\beta$ with $\beta'$, Fig.~\ref{fig:BMO}. 

\phantom{x}

By this operation, the new 
banded loop system $(\mbox{\boldmath $\kappa$}, {\mbox{\boldmath $\beta$}}')$ 
is a renewal embedding of the original banded loop system 
$(\mbox{\boldmath $\kappa$}, {\mbox{\boldmath $\beta$}})$ and 
has as the surgery link a new union $k'_0$ of the same links $k$ and 
${\mathbf o}$, not necessarily the split sum $k+{\mathbf o}$,  
because the band system ${\mbox{\boldmath $\beta$}}'$ is isotopic to 
${\mbox{\boldmath $\beta$}}$ if ${\mbox{\boldmath $\alpha$}}^c_1$ is forgotten. 
In the final stage of this paper, the surgery link $k'_0$  will 
have $k\cap {\mathbf d}=\emptyset$, so that $k'_0$ will be 
the split sum $k+{\mathbf o}$, because ${\mathbf o}=\partial{\mathbf d}$. 

To achieve a situation where the Band Move Operation can be applied, 
the follwing concept is needed. 
A {\it splitting} of a banded loop system 
$(\mbox{\boldmath $\kappa$}, {\mbox{\boldmath $\beta$}})$ 
is a banded loop system $(\mbox{\boldmath $\kappa$}^*, {\mbox{\boldmath $\beta$}}^*)$ 
such that a based disk system $\mbox{\boldmath $\delta$}^*$ for 
$\mbox{\boldmath $\kappa$}^*$
is obtained from a based disk system $\mbox{\boldmath $\delta$}$ for 
$\mbox{\boldmath $\kappa$}$ by splitting along a disjoint proper arc system 
${\mbox{\boldmath $\gamma$}}$ in $\mbox{\boldmath $\delta$}$ not meeting ${\mathbf o}$ 
and ${\mbox{\boldmath $\beta$}}$,  
and the band system ${\mbox{\boldmath $\beta$}}^*$ 
is obtained from the band system ${\mbox{\boldmath $\beta$}}$ by adding 
the band system ${\mbox{\boldmath $\beta_{\gamma}$}}$ thickening 
${\mbox{\boldmath $\gamma$}}$. This splitting operation comes from 
Fission-Fusion move of a banded loop system, \cite{Ka1}. 
After some splittings of a banded loop system,  a situation where the Band Move Operation 
can be applied is realized by a replacement of the based disk system and an isotopic deformation of the band system.

The following assertion is used.

\phantom{x}

\noindent{\bf (2.3)} If there is a splitting 
$(\mbox{\boldmath $\kappa$}^*, {\mbox{\boldmath $\beta$}}^*)$ of 
a banded loop system 
$(\mbox{\boldmath $\kappa$}, {\mbox{\boldmath $\beta$}})$ 
with surgery knot $k_0$ a union of $k$ and ${\mathbf o}$
such that $\mbox{\boldmath $\kappa$}^*$ does not meet the interior of the extra disk system 
${\mathbf d}$, then 
there is a renewal embedding $(\mbox{\boldmath $\kappa$}', {\mbox{\boldmath $\beta$}}')$ 
of $(\mbox{\boldmath $\kappa$}, {\mbox{\boldmath $\beta$}})$ 
such that $(\mbox{\boldmath $\kappa$}', {\mbox{\boldmath $\beta$}}')$ 
does not meet the interior of ${\mathbf d}$ and has the surgery knot 
$k'_0=k+{\mathbf o}$.

\phantom{x}

\noindent{\it Proof of (2.3).} Since ${\mbox{\boldmath $\kappa$}}^*$ does not meet 
the interior of ${\mathbf d}$, there is a based disk system 
${\mbox{\boldmath $\delta$}}^*$ 
for ${\mbox{\boldmath $\kappa$}}^*$ not meeting the interior of ${\mathbf d}$.
The band system ${\mbox{\boldmath $\beta$}}^*$ 
transversely meets the interior of ${\mathbf d}$ with transverse arc system $A$. 
Let  ${\mbox{\boldmath $\delta$}}^*_1$ be the sub-system 
of ${\mbox{\boldmath $\delta$}}^*$ containing the complementary arc system 
${\mbox{\boldmath $\alpha$}}^c_1$ in the boundary, and 
$N({\mbox{\boldmath $\alpha$}}^c_1)$ a boundary collar disk system of  
${\mbox{\boldmath $\alpha$}}^c_1$ in ${\mbox{\boldmath $\delta$}}^*_1$. 
The Band Move Operation means that the band system ${\mbox{\boldmath $\beta$}}^*$ is 
deformed  
so that the transverse arc system $A$ moves from the interior of ${\mathbf d}$ 
into the interior of $N({\mbox{\boldmath $\alpha$}}^c_1)$. 
Then by changing  the band
system ${\mbox{\boldmath $\beta_{\gamma}$}}$  back into the arc system 
${\mbox{\boldmath $\gamma$}}$, 
the banded loop system 
$(\mbox{\boldmath $\kappa$}^*, {\mbox{\boldmath $\beta$}}^*)$ is  changed back to a pair 
$(\mbox{\boldmath $\kappa$}', {\mbox{\boldmath $\beta$}}')$, where 
the loop system $\mbox{\boldmath $\kappa$}'$ bounds an immersed disk system 
$\mbox{\boldmath $\delta$}'$ obtained from the based disk system 
${\mbox{\boldmath $\delta$}}$  by moving a transverse arc system of 
${\mbox{\boldmath $\beta_{\gamma}$}}$ into the interior of 
$N({\mbox{\boldmath $\alpha$}}^c_1)$. 
The immersed disk system $\mbox{\boldmath $\delta$}'$ is deformed into a  disjoint disk system by repeatedly pulling the band in ${\mbox{\boldmath $\beta_{\gamma}$}}$  connecting to an outer most disk of $\mbox{\boldmath $\delta$}^*$ or passing 
the outer most disk of $\mbox{\boldmath $\delta$}^*$ 
through $N({\mbox{\boldmath $\alpha$}}^c_1)$ in order to eliminate the nearest 
transverse arc  of the band.
This means that the  loop system 
$\mbox{\boldmath $\kappa$}'$ is  a trivial link and   
$(\mbox{\boldmath $\kappa$}', {\mbox{\boldmath $\beta$}}')$ is a banded loop system. 
Thus,  there is a renewal embedding 
$(\mbox{\boldmath $\kappa$}', {\mbox{\boldmath $\beta$}}')$ 
of $(\mbox{\boldmath $\kappa$}, {\mbox{\boldmath $\beta$}})$ 
which does not meet the interior of ${\mathbf d}$. The surgery knot $k'_0$ is necessarily 
the split sum $k+{\mathbf o}$ since $\partial {\mathbf d}={\mathbf o}$.
This completes the proof of (2.3).

\phantom{x}

By using (2.2) and (2.3), the following assertion is shown.

\phantom{x}

\noindent{\bf (2.4)} There is a renewal embedding 
$(\mbox{\boldmath $\kappa$}', {\mbox{\boldmath $\beta$}}')$ of 
every banded loop system 
$(\mbox{\boldmath $\kappa$}, {\mbox{\boldmath $\beta$}})$ in 
${\mathbf R}^3$ with surgery link $k_0=k+{\mathbf o}$ 
such that  $(\mbox{\boldmath $\kappa$}', {\mbox{\boldmath $\beta$}}')$ 
does not meet the interior of ${\mathbf d}$ and has the surgery knot 
$k'_0=k+{\mathbf o}$.

\phantom{x}

\noindent{\it Proof of (2.4).} 
By (2.2), a based disk system ${\mbox{\boldmath $\delta$}}$ of 
$\mbox{\boldmath $\kappa$}$ transversely 
meets the extra disk system ${\mathbf d}$ with interior simple arcs or  normal proper arcs in ${\mathbf d}$ except for the complementary arc system 
${\mbox{\boldmath $\alpha$}}_1^c$. 
Let $A$ be the interior arc system which is made disjoint from 
${\mbox{\boldmath $\beta$}}$ 
by isotopic deformations of ${\mbox{\boldmath $\beta$}}$ respecting  the arc system 
${\mbox{\boldmath $\alpha$}}_1$ and the loop system 
$\mbox{\boldmath $\kappa$}$.
By taking a splitting of $(\mbox{\boldmath $\kappa$}, {\mbox{\boldmath $\beta$}})$ 
along $A$, it is considered that 
the based disk system ${\mbox{\boldmath $\delta$}}$ transversely 
meets ${\mathbf d}$ only with  normal proper arcs in ${\mathbf d}$ except for 
${\mbox{\boldmath $\alpha$}}_1^c$. 
Then $\mbox{\boldmath $\kappa$}$ does not meet the interior of the extra disk system 
${\mathbf d}$. By (2.3), the proof of (2.4) is completed.

\phantom{x}

Let 
$(\mbox{\boldmath $\kappa$}, {\mbox{\boldmath $\beta$}})$ be a  banded loop system 
a banded loop system with surgery link 
$k_0=k+{\mathbf o}$ such that $P(F^1_0)$ is equivalent to $F$. 
By (2.4), there is a renewal embedding 
$(\mbox{\boldmath $\kappa$}', {\mbox{\boldmath $\beta$}}')$ 
such that $(\mbox{\boldmath $\kappa$}', {\mbox{\boldmath $\beta$}}')$ does 
not meet the interior of the extra disk system ${\mathbf d}$, and  
has the surgery link $k+{\mathbf o}$. 
Let ${\mathbf b}'$ be the band system  dual to the band system 
${\mbox{\boldmath $\beta$}}'$, and 
${\mathbf b}'_1$ the band sub-system of ${\mathbf b}'$ such that ${\mathbf b}'_1$ 
connects to $ {\mathbf o}$ with just one band for  every component of $ {\mathbf o}$. 
Let ${\mathbf b}'_2={\mathbf b}'\setminus {\mathbf b}'_1$. 
Since ${\mathbf b}'_1$ does not meet the interior of ${\mathbf d}$, 
the surgery link of the banded link $(k+{\mathbf o}, {\mathbf b}'_1)$ 
is equivalent to the link $k$ and the upper-closed realizing surface 
$\mbox{ucl}(F^1_0)'$ of the banded link $(k, {\mathbf b}'_2)$ 
is equivalent to the proper realizing surface $P(F_0^1)'$ of 
$(\mbox{\boldmath $\kappa$}', {\mbox{\boldmath $\beta$}}')$
which is a ribbon surface in ${\mathbf R}^4_+$ and is a renewal embedding of 
the proper realizing surface $P(F_0^1)$ of the  banded loop system 
$(\mbox{\boldmath $\kappa$}, {\mbox{\boldmath $\beta$}})$ with 
the surgery link $k+{\mathbf o}$. 
Since $P(F_0^1)$ is equivalent to $F$ in ${\mathbf R}^4_+$ 
and $\mbox{ucl}(F^1_0)'$ is a ribbon surface with $\partial\mbox{ucl}(F^1_0)'=\partial F=k$,
there is a boundary-relative renewal embedding from $\mbox{ucl}(F^1_0)'$ to $F$.
This completes the proof of the classical ribbon theorem. 

\phantom{x}

\noindent{\bf Acknowledgements.} The author thanks to a referee for suggesting 
the content for making easier content.
This paper is motivated by T. Shibuya's comments pointing out insufficient explanation on 
Lemma 2.3 in \cite{Ka2} (Corollary~2 in this paper). 
This work was partly supported by JSPS KAKENHI Grant Number JP21H00978 and MEXT 
Promotion of Distinctive Joint Research Center Program JPMXP0723833165.

\phantom{x}

\end{document}